\magnification=\magstephalf
\input eplain
\input BMmacs

\phantom{x}\vskip .5cm
\centerline{{\titlefont A REMARK ON THE TATE CONJECTURE}}
\vskip .5cm

\centerline{{\byfont by}\quad {\namefont Ben Moonen}}
\vskip 1cm

{\eightpoint
\noindent
{\bf Abstract.} The strong version of the Tate conjecture has two parts: an assertion~(S) about semisimplicity of Galois representations, and an assertion~(T) which says that every Tate class is algebraic. We show that in characteristic~$0$, (T) implies~(S). In characteristic~$p$ an analogous result is true under stronger assumptions.
\medskip

\noindent
{\it AMS 2010 Mathematics subject classification:\/} 14C15, 14C25, 14F20\par}
\vskip 8mm

\noindent
Let $K$ be a finitely generated field of characteristic~$0$, and fix an arbitrary prime number~$\ell$. Choose an algebraic closure $K \subset \Kbar$ and let $\Gamma_K = \Gal(\Kbar/K)$.

Let $X$ be a smooth projective $K$-scheme. Given integers $i\geq 0$ and~$n$, we write $H^i(X)(n)$ for the $\ell$-adic \'etale cohomology $H^i(X_\Kbar,\Ql(n))$, which comes equipped with a continuous action of~$\Gamma_K$. A class $\xi \in H^i(X)(n)$ is called a Tate class if its stabilizer in~$\Gamma_K$ is an open subgroup. If $\xi \neq 0$ is a Tate class then necessarily $i=2n$. Let $\Tate^n(X) \subset H^{2n}(X)(n)$ denote the subspace of Tate classes.

We are interested in the following two conjectures:
\item{(S)} For every $X/K$, $i$ and~$n$ as above, $H^i(X)(n)$ is a semisimple representation of~$\Gamma_K$.
\item{(T)} For every $X/K$ as above and $n \geq 0$, the cycle class map $\CH^n(X_\Kbar) \otimes \Ql \to \Tate^n(X)$ is surjective.
\smallskip

The conjunction of these two conjectures is the (strong version of the) Tate conjecture, over fields of characteristic~$0$. The semisimplicity of $H^i(X)(n)$ as a Galois representation is equivalent to the semisimplicity of $H^i(X) = H^i(X)(0)$; hence in Conjecture~(S) we may restrict our attention to the case $n=0$.

The goal of this note is to prove the following result.
\medskip

\noindent
{\it Theorem~1. --- Assume Conjecture~{\rm (T)} is true. Then also Conjecture~{\rm (S)} is true.\/}
\medskip

In what follows we assume Conjecture~(T) is true. We take $X/K$ and $i\geq 0$ as above, and we abbreviate $H^i(X)$ to~$H$. Our goal is to prove that the Galois representation $\rho \colon \Gamma_K \to \GL(H)$ is semisimple. Let $G \subset \GL(H)$ be the Zariski closure of the image of~$\rho$; then the semisimplicity of~$\rho$ is equivalent to the assertion that the identity component~$G^0$ is a reductive group. If we replace $K$ by a finitely generated field extension, the group~$G^0$ does not change; in particular, $\rho$ is semisimple if and only if its restriction to an open subgroup of~$\Gamma_K$ is semisimple.
\medskip

{\it Reduction to the case $K= \mQ$.} First we can reduce to the case that $K$ is a number field. For this we use a (standard) specialization argument due to Serre; see the first section of~[\ref{Serre2Ribet}] or [\ref{SerreMW}], Section~10.6. Next, let $X^\prime$ be $X$ viewed as a scheme over~$\mQ$ (so: $X \to \Spec(K) \to \Spec(\mQ)$). Then $H^\prime = H^i(X^\prime_\Qbar,\Ql)$ is the Galois representation obtained from~$H$ by induction from $\Gamma_K$ to~$\Gamma_\mQ$. If $H^\prime$ is semisimple as a representation of~$\Gamma_\mQ$, so is its restriction to~$\Gamma_K$, and it follows that $H$ is a semisimple representation of~$\Gamma_K$, too.
\medskip

We now assume that $X$ is smooth projective over~$\mQ$. Write $\Gamma = \Gamma_\mQ$. Let $V \subset H = H^i(X)$ be a $\Gamma$-submodule. We are done if we can show that $V$ has a complement in~$H$ that is stable under the action of an open subgroup of~$\Gamma$. Let $m = \dim(V)$, and write $Y = X^m$. Then $H^{mi}(Y_\Qbar,\Ql)$, as a representation of~$\Gamma$, contains a copy of $\Alt^m(H) \cong \wedge^m H$, which in turn contains the line~$\Alt^m(V)$.

As $V$ contains a $\Gamma$-stable $\Zl$-lattice, the Galois action on $\Alt^m(V)$ is given by a character $\psi \colon \Gamma \to \Zl^\times$. We claim that the restriction of~$\psi$ to an open subgroup of~$\Gamma$ is equal to an integral power of the $\ell$-adic cyclotomic character~$\chi$. As this is trivially true if $\psi$ has finite image, we may assume $\psi(\Gamma)$ to be infinite. Let $\mQ \subset L$ be the field extension corresponding to $\Ker(\psi)$. Then $L$ contains the cyclotomic $\Zl$-extension $\mQ_\infty$ of~$\mQ$ as a subfield of finite index. As $\ell$ is totally ramified in~$\mQ_\infty$, it follows that the inertia group $I_\lambda \subset \Gal(L/\mQ)$ of a prime~$\lambda$ above~$\ell$ is an open subgroup. On the other hand, as proven by Faltings~[\ref{Faltings}], $H$ is a Hodge--Tate representation of~$\Gamma_{\Ql}$. This implies that~$V$, and then also $\Alt^m(V)$ is Hodge--Tate. As $\Alt^m(V)$ is $1$-dimensional, it follows (see for instance [\ref{Fontaine}], Remark~3.9(\romannumeral4)) that there is an open subgroup of~$I_\lambda$ on which~$\psi$ is given by an integral power of~$\chi$. This gives our claim.

To conclude the argument we now consider, for $F$ a field of characteristic~$0$, the category $\Mot(F;\Ql)$ of motives over~$F$ in the sense of Andr\'e~[\ref{YA}], with coefficients in~$\Ql$. This is a semisimple $\Ql$-linear neutral Tannakian category. Let $\cG_{\mot,F}$ denote the motivic Galois group corresponding with the fibre functor given by $\ell$-adic cohomology, which is a (non-connected) pro-reductive group. The category $\Mot(F;\Ql)$ is equivalent with the category $\Rep(\cG_{\mot,F};\Ql)$ of representations of~$\cG_{\mot,F}$ on finite-dimensional $\Ql$-vector spaces. With $X/\mQ$ as above, $H = H^i(X)$ is (the $\ell$-adic realization of) a submotive of~$X$ over~$\mQ$. If $F$ is a number field, we therefore have a natural action of $\cG_{\mot,F}$ on~$H$.

As the Galois representation $H^{mi}(Y_\Qbar,\Ql)$ is pure of weight~$mi$, it follows from what we have shown that $mi$ is even, say $mi=2k$, and that there exists a number field~$F$ such that $\Alt^m(V) \cong \Ql(-k)$ as representations of~$\Gamma_F$. Conjecture~(T) then implies that the line $\Alt^m(V)(k) \subset H^{mi}(Y_\Qbar,\Ql)(k) = H^{2k}(Y_\Qbar,\Ql)(k)$ is spanned by the cohomology class of an algebraic cycle in $\CH^k(Y_F)$. In particular, $\Alt^m(V)$ is a motivated subspace of $\Alt^m(H)$ over~$F$, i.e., $\Alt^m(V)$ is a $\cG_{\mot,F}$-submodule of~$\Alt^m(H)$. This implies that $V$ is a $\cG_{\mot,F}$-submodule of~$H$. (The map $V \mapsto \Alt^m(V)$ gives the Pl\"ucker embedding of the Grassmannian of $m$-planes in~$H$.) As $\Mot(F;\Ql)$ is a semisimple category, it follows that $H$ decomposes as $V \oplus V^\prime$ as a motive over~$F$, and hence also as a representation of~$\Gamma_F$. \QED
\medskip

\noindent
{\it Remarks.\/} --- (\romannumeral1) Let $X$ be projective smooth over a number field~$K$ that is Galois over~$\mQ$. Consider varieties~$Y$ obtained as products of finitely many conjugates of~$X$ over~$\mQ$. The proof shows that if all $\Gamma_K$-invariant classes in all $H(Y)\bigl(n\bigr)$ are $\Ql$-linear combinations of motivated cycles (in the sense of Andr\'e~[\ref{YA}]) then $H(X)$ is a semisimple $\Gamma_K$-representation. Instead of using motivated cycles, we could also work with absolute Hodge classes, as introduced by Deligne.

(\romannumeral2) The reduction to $\mQ$ as a base field is essential. In general, if $X/K$ is smooth projective and $V \subset H^i(X)$ is an $m$-dimensional $\Gamma_K$-submodule, it is not true that $\Alt^m(V)$ is spanned by a Tate class. For instance, suppose $X$ is an elliptic curve whose endomorphism algebra is an imaginary quadratic field~$k$, and suppose $\ell$ splits in~$k$; then $H^1(X_\Kbar,\Ql)$ is the direct sum of two $1$-dimensional Galois representations that are not (potentially) Tate twists. Note that this case gives rise to $\Zl$-extensions of~$K$ in which the primes above~$\ell$ have finite ramification.

(\romannumeral3) Assuming (T) is true, Theorem~1 implies that the image of the Galois representation $\rho \colon \Gamma_K \to \GL(H)$ is Zariski-dense in the motivic Galois group of the motive~$H^i(X)$.

(\romannumeral4) As Peter Scholze pointed out to us, the proof may be concluded without reference to motives, as follows. As (T) implies that the line $\Alt^m(V)(k)$ is generated by a Hodge class in the Betti cohomology, the Hodge-Riemann relations together with hard Lefschetz imply that there exists a $\Gamma_F$-equivariant splitting $s \colon \Alt^m(H) \to \Alt^m(V) \cong \Ql(-k)$. Via the composition
$$
V^\vee(-k) \otimes H \mapright{\sim} \Alt^{m-1}(V) \otimes H \longrightarrow \Alt^{m-1}(H) \otimes H \longrightarrow \Alt^m(H) \mapright{s} \Alt^m(V)
$$
we get a $\Gamma_F$-equivariant map $V^\vee \otimes H \to \Ql$; whence a splitting $H \to V$ of the inclusion $V \hookrightarrow H$.
\bigskip

To extend the argument in the proof to fields of characteristic~$p$ seems to require further assumptions, e.g., that homological and numerical equivalence agree. But as we learned from Marco D'Addezio, under such assumptions, (S) can be deduced from results in the literature:
\medskip

\noindent
{\it Theorem~2. --- Let $p$ and $\ell$ be prime numbers with $\ell \neq p$. Assume Conjecture~{\rm (T)} is true over finite fields of characteristic~$p$ and that ${\equiv_\num} = {\sim_\hom}$ over such fields. Then Conjecture~{\rm (S)} is true in characteristic~$p$, i.e., for $X$ smooth projective over a finitely generated field~$K$ of characteristic~$p$, the $\Gamma_K$-representation $H^*(X_\Kbar,\Ql)$ is semisimple.\/}
\medskip

\noindent
{\it Proof.} Over a finite field, this follows from [\ref{Milne}], Prop.~8.4 and Remark~8.6. The case of a finitely generated field then follows from the main result of~[\ref{Fu}]. \QED
\medskip

Over a finite field, assuming ${\equiv_\num} = {\sim_\hom}$, [\ref{Milne}] in fact shows that the Galois action on $H^*(X_\Kbar,\Ql)$ is semisimple if the Tate conjecture is true for cycles of codimension~$\dim(X)$ on~$X^2$.
\bigskip

\noindent
{\it Acknowledgements.} --- I thank Marco D'Addezio for his comments and suggestions with regard to the situation in positive characteristic, and in particular for the references to~[\ref{Fu}] and~[\ref{Milne}]. I also thank Anna Cadoret for her remarks on a first version of this paper, and Peter Scholze for Remark~(\romannumeral4).

\bigskip

{\eightpoint
\goodbreak\centerline{{\bf References}}%
\nobreak\vskip.5\bigskipamount plus 2pt minus 1pt%

\bibitem{YA}
Y.~Andr\'e,
{\it Pour une th\'eorie inconditionnelle des motifs.}
Inst.\ Hautes \'Etudes Sci. Publ. Math. 83 (1996), 5--49.

\bibitem{Faltings}
G.~Faltings,
{\it $p$-adic Hodge theory.}
J.\ Amer.\ Math.\ Soc.\ 1 (1988), no.~1, 255--299.

\bibitem{Fontaine}
J.-M.~Fontaine,
{\it Repr\'esentations $p$-adiques semi-stables.}
In: P\'eriodes $p$-adiques (Bures-sur-Yvette, 1988).
Ast\'erisque 223 (1994), 113--184.

\bibitem{Fu}
L.~Fu,
{\it On the semisimplicity of pure sheaves.}
Proc.\ Amer.\ Math.\ Soc.\ 127 (1999), no.~9, 2529--2533.

\bibitem{Milne}
J.~Milne,
{\it Values of zeta functions of varieties over finite fields.}
Amer.\ J.\ Math.\ 108 (1986), no.~2, 297--360.

\bibitem{Serre2Ribet}
J-P.~Serre,
{\it Lettres \`a Ken Ribet du 1/1/1981 et du 29/1/1981.}
In: \OE uvres, Volume~IV, number 133.
Springer-Verlag, Berlin, 2000.

\bibitem{SerreMW}
J-P.~Serre,
{\it Lectures on the Mordell-Weil theorem.}
Aspects of Math., Vol.~E15.
Friedr.\ Vieweg \& Sohn, Braunschweig, 1997.

\bigskip

\noindent
Radboud University, IMAPP, PO Box 9010, 6500GL Nijmegen, The Netherlands.

\noindent
{\tt b.moonen@science.ru.nl}
}

\bye